\magnification=1200

{\centerline {\bf Degree Estimates for Polynomials Constant on a Hyperplane}}
\bigskip

John P. D'Angelo, Ji{\v r}\'\i\ Lebl, and Han Peters.

Addresses:

D'Angelo (jpda@math.uiuc.edu), Dept. of Mathematics,
Univ. of Illinois, Urbana IL 61801.

Lebl (jlebl@math.ucsd.edu), Dept. of Mathematics, Univ. of California at 
San Diego, La Jolla CA 92093.

Peters (peters@math.wisc.edu), Dept. of Mathematics, 
Univ. of Wisconsin, Madison WI 53706.

\bigskip

{\bf I. Introduction}
\medskip

We are interested in the complexity of real-valued polynomials, 
defined on real Euclidean space
${\bf R}^n$, that are constant 
on a hyperplane. This issue arises as a simplified version
of a difficult question in CR Geometry, which we discuss briefly 
below and in Section VI.
We intend to fully address the CR issues in a subsequent paper.

Let $H$ denote the hyperplane in ${\bf R}^n$
defined by $\{x: s(x) = \sum_{j=1}^n x_j = 1\}$. 
We write ${\bf R}[x] = {\bf R}[x_1,...,x_n]$ for the ring of
real-valued polynomials in $n$ real variables.
Suppose $p \in {\bf R}[x]$ and that $p$ is constant on $H$.
How complicated can $p$ be?
Two possible measurements of the complexity of a polynomial 
are its degree $d$ and the number $N$
of its distinct monomials. We always have the standard estimate

$$ N \le {n+d \choose n}, \eqno (1) $$
which estimates $d$ from below. Even when $p$ is constant on $H$,
no upper estimate for $d$ in terms of $N$
is possible without additional assumptions.
For example, for $d\ge 2$, consider
$$ p(x) = x_1^{d-1} s(x) - x_1^{d-1} + 1. \eqno (2) $$
It is evident that $p=1$ on $H$, that $p$ has $n+2$ distinct monomials, 
and that its degree $d$ can be arbitrarily large. On the other hand, such degree
estimates become possible when we assume that $n\ge 2$ and that
the coefficients of $p$ are nonnegative. We prove such results in this paper.

Before describing our results
we briefly discuss the motivation behind them. See Section VI for additional information.
In a future paper we will say more
about this connection with CR geometry.
Let $f:{\bf C}^n \to {\bf C}^N$ be a rational mapping such that $f$ maps the unit ball
in its domain properly to the unit ball in its target.
It follows that $f$ maps the unit sphere in ${\bf C}^n$ to the unit sphere
in ${\bf C}^N$. For $n\ge2$, the work of
Forstneric [F1] implies that the degree of $f$ is bounded in terms of $n$ and $N$. 
The bound in [F1] is not sharp, and finding a sharp bound seems to be difficult.
Meylan [M] has improved the bound when $n=2$.

The problem simplifies somewhat by assuming that $f$ is a monomial mapping; that is,
$f$ is a polynomial mapping for which (after a coordinate change if necessary)
each component is a monomial. The condition $||f(z)||^2=1$ on $||z||^2=1$ then 
depends upon only the real variables $|z_1|^2,..., |z_n|^2$, and all coefficients 
involved appear as $|c|^2$ for complex numbers $c$. The relationship between the 
degree of $f$ and the domain and target dimensions then becomes the 
combinatorial issue described in Problem 1 below.

We need to consider various subsets of ${\bf R}[x_1,...,x_n]$.
Let ${\cal J}(n)$ denote the subset of polynomials $p$ in ${\bf R}[x_1,...,x_n]$
for which $p(x)=1$ on the hyperplane $H$. The set ${\cal J}(n)$  
is closed under multiplication,
convex combinations, and the operation $X$ described in Section II. Let
${\cal P}(n)$ denote those polynomials in ${\bf R}[x_1,...,x_n]$
whose coefficients are {\it nonnegative}. The set ${\cal P}(n)$ is closed under
addition and multiplication.
Let ${\cal P}(n,d)$ denote the subset of ${\cal P}(n)$ 
whose elements are of degree $d$. The crucial sets for us are
${\cal H}(n)$ and ${\cal H}(n,d)$: 

$${\cal H}(n)= {\cal J}(n) \cap {\cal P}(n)$$ 
$${\cal H}(n,d)= {\cal J}(n) \cap {\cal P}(n,d). $$

Thus the elements of ${\cal H}(n,d)$ are polynomials of degree $d$ in $n$ real
variables, with nonnegative coefficients, and whose values are $1$ on the set
$\sum x_j =1$. 
For $p \in {\bf R}[x]$, we write $N(p)$ for the number of distinct
monomials occurring in $p$. Our goal is to prove sharp
estimates relating the degree of $p$ with $N(p)$ when $p\in {\cal H}(n)$.

\medskip {\bf Problem 1}. Assume $n\ge 2$.
For $p \in {\cal H}(n)$, find a sharp upper bound for $d(p)$
in terms of $N(p)$ and $n$.
\medskip

There is no such upper bound when $n=1$, as we note in Section II. When
$n=2$, the sharp upper bound is given by $d(p) \le 2 N(p) -3$, a result from [DKR]
we discuss also in Section II. For $n\ge 3$ the first author has conjectured the bound

$$ d(p) \le  { N(p)-1 \over n-1}. \eqno (3) $$
Example 4 provides polynomials of each degree
where equality holds in (3). 

In Proposition 4 from Section III we pullback
to the two-dimensional case via a Veronese mapping
to obtain a general but crude bound.
For $n\ge 2$ and $p \in {\cal H}(n,d)$ we obtain 

$$ d(p) \le {2N(p)-3 \over n-1}. \eqno (4) $$

This result is not sharp unless $n=2$. In Section IV we improve (4)
by pulling back via the optimal mappings in two dimensions.
In Theorem 1 we obtain 

$$ d(p) \le {2n(2N(p)-3) \over 3n^2-3n-2} \le {4 \over 3} {2N(p)-3 \over 2n-3}. 
\eqno (5) $$

In Theorem 2 of Section V we prove our main result:
for $n$ sufficiently large compared with $d$, the estimate (3) holds, 
and we find all polynomials
for which equality holds in (3). We remark now, 
and demonstrate later, that when $n=3$ for
example, there are additional polynomials for which equality holds. 
It is therefore reasonable to think of Theorem 2 as a stabilization result; certain
complicated issues arise in low dimensions, but become irrelevant as the dimension
$n$ rises. In Corollary 2 of Section IV we also lend support to the conjecture. When $n\ge3$
we show that the conjecture holds for degree up to $4$. We show also there that the conjecture holds
when $N < 4n-3$.

We summarize our work. In Theorem 1 we prove a general bound 
which is not sharp unless $n=2$.
Lemmas 4 and 5 show how to sharpen that bound in specific situations.
In Corollary 2 we prove a sharp bound for all $n$ when either $d\le 4$ or $N < 4n-3$.
In Theorem 2 we establish the sharp bound when $n$ is sufficiently large given $d$.

We close the introduction with one additional comment.
When $p \in {\cal J}(n)$, the function
$Q(p)$, defined by
$$ Q(p) = {p-1 \over s-1} \eqno (6) $$
is a polynomial; 
its coefficients need not be nonnegative even if $p \in {\cal H}(n)$. 
The polynomial $Q$ plays a crucial role in the proof
in two dimensions, and it therefore plays an implicit role here.
Perhaps some of our results can be better understood in terms of $Q(p)$.

The first author posed Problem 1 at the workshop on CR Geometry held at MSRI
in July 2005; the other two authors attended that workshop and began working on
it at that time. All three authors acknowledge MSRI.  
The three authors obtained one of the results here and
put the finishing touches on this paper at the workshop on CR Geometry at AIM, 
September 2006.  All three authors thus acknowledge AIM. The first author also
acknowledges NSF Grant DMS 0500765.

\medskip
{\bf II. The situations in one and two dimensions}
\medskip

The situation in one dimension is not interesting, so we dispense with it now, 
and assume thereafter that $n\ge2$. When $n=1$, we note that $p\in {\cal H}(1)$ 
when $p$ has nonnegative coefficients and $p(1)=1$. 
The particular polynomial $p(x_1) = x_1^d$ lies in ${\cal H}(1,d)$ 
and $N(p)=1$. Furthermore, for any fixed value of $N$, we can find a polynomial $p$ 
of arbitrarily large degree with $N(p)=N$. Thus no upper bound for $d(p)$ is possible.

When $n=2$ a sharp result is known [DKR].
\medskip

{\bf Theorem 0}. Let $p$ be a polynomial in two real variables $(x,y)$
such that

1) $p(x,y)= 1$ when $x+y=1$, and

2) each coefficient of $p$ is nonnegative.

\noindent 
Let $N$ be the number of distinct monomials in $p$, and let $d$ be the degree
of $p$. Then $d \le 2N-3$.  Furthermore, for each $N\ge 2$,
there is a polynomial $p_d$ satisfying 1) and 2) whose degree is $2N-3$.
\medskip
The estimate $d \le 2N-3$ can of course be rewritten
$N \ge {d+3 \over 2}$. The proof of Theorem 0 shows that a
slightly stronger conclusion holds.
If $p$ satisfies 1) and 2) then $p$ must have
at least ${d-1 \over 2}$ mixed terms (those containing both $x$ and $y$)
and at least two pure terms.

\medskip
There is an interesting family
of polynomials providing the sharp bound in Theorem 0.
The polynomials in this family have integer coefficients, they are group-invariant, and
they exhibit many interesting combinatorial and number-theoretic properties. We
mention for example that $p_d(x,y) \cong x^d + y^d$ if and only if $d$ is prime.
See [D1,D2, D3, D5, DKR] for this fact and much additional information. 
Here is an explicit formula for these polynomials for $d$ odd:

$$ p_d(x,y) = y^d  + ({{ x + {\sqrt {x^2 + 4 y}}} \over 2})^d + 
({{ x - {\sqrt {x^2 + 4 y}}} \over 2})^d.  \eqno (7) $$

We also provide a recurrence formula relating these polynomials as the degree varies.
Put $g_0(x,y) = x$ and $g_1(x,y) = x^3 + 3xy$.
Define $g_{k+2}$ and then $p_{2k+1}$ by
$$ g_{k+2} (x,y) = (x^2 + 2y)g_{k+1}(x,y) - y^2 g_k(x,y)$$
$$p_{2k+1}(x,y) = g_k(x,y) + y^{2k+1}. \eqno (8) $$

The equations in (8) determine the polynomials in (7).
For odd $d$ the polynomial defined by (7) has precisely ${d+3 \over 2}$ terms,
and thus the bound in Theorem 0 is sharp.
We can obtain a second sharp example by interchanging
the roles of $x$ and $y$. Other examples exhibiting the sharp bound exist for
some but not all $N$. See Example 3 where $N=5$.

Each $p_{2r+1}$ is group-invariant; we have
$ p_{2r+1} (\omega x, \omega^2 y) = p_{2r+1} (x,y)$ 
whenever $\omega$ is a $2r+1$-st root of unity.
There are analogous group-invariant polynomials for even degree, but
these have a single negative
coefficient, and we will not discuss them in this paper.

\medskip

The proof of the inequality  $d \le 2N-3$ from Theorem 0 
is quite complicated. It relies on
an analysis of certain directed graphs arising from the Newton diagram of the
polynomial $Q(p)$ and their interaction with Proposition 1 below.

We close this section
by indicating how one can use Theorem 0 to study the higher-dimensional case.
Let $\phi : {\bf R}^2 \to {\bf R}^n$ be a polynomial mapping, and suppose
that $\phi$ maps the line defined by $u+v=1$ to the hyperplane $H$. If
$p \in {\cal J}(n)$, then the composite map $\phi^*(p)$ is in  ${\cal J}(2)$.
To see this fact, observe for $u+v=1$ that 

$$ \phi^*(p) (u,v) = p(\phi(u,v)) = 1 \eqno (9) $$
because $p=1$ on $H$.

We will apply this idea of pulling back to two dimensions
for various functions $\phi$. We give some examples.
Assume $n \ge 3$. For $i \ne j$,
set $x_i = u$, set $x_j = v$, and set
$x_k = 0$ otherwise. Another possibility is to set $k$ of the variables equal 
to ${u \over k}$, set $l$ of the other variables equal to ${v \over l}$, 
and set the remaining variables equal to zero. In these cases $\phi$ is linear. 
In the proof of Proposition 4 from Section III we let $\phi$ be a Veronese mapping;
in that proof $\phi$ is homogeneous of degree larger than one. One can also
gain information by pulling back via more complicated mappings. See 
Sections IV and V  for details.

\medskip
{\bf III. General Information}
\medskip

We begin with several formal algebraic observations. Suppose that $p \in {\cal J}$, 
and that $u$ is an arbitrary polynomial. We define a polynomial $X_u (p)$
by

$$ X_u(p) = p-u + su. \eqno (10) $$
When $p \in {\cal J}$ we can always write $p = (1 - Q) + sQ$ where $Q$ is as in
(6), and thus $p=X_Q(1)$.  In general we will drop the dependence on $u$ from the notation and write $X(p)$ for $X_u(p)$.  The following simple but crucial result
suggests decomposing elements in ${\cal H}$ using the operation in (10).

\medskip
{\bf Lemma 1}. Suppose $p \in {\cal J}$ and $u$ is an arbitrary polynomial.
Define $X(p)$ by (10). Then $X(p) \in {\cal J}$. Suppose $p \in {\cal H}$ 
and also suppose that both $u$ and $p-u$ are in ${\cal P}$. Then $X(p) \in {\cal H}$.

Proof. It is immediate from (10) that $X(p) = p-u +u = p$ on the set $s=1$, and hence $X(p)  \in {\cal J}$.
Suppose that both $u$ and $p-u$ are in ${\cal P}$. Also $s \in {\cal P}$. Since
${\cal P}$ is closed under addition and
multiplication, it follows that $X(p) \in {\cal P}$. Since we have shown that $X(p) \in {\cal J}$ as well,
$X(p) \in {\cal H}$. $\spadesuit$

\medskip
Our concern with nonnegative coefficients leads us to make the following definition.

\medskip
{\bf Definition 1}. Suppose that $p,g \in {\cal P}(n)$. We say that $g \subset p$
if $p-g \in {\cal P}(n)$. In other words $g \subset p$ holds if and only if
both $g$ and $p-g$ have nonnegative coefficients. 
We call $g$ a {\it subpolynomial} of $p$.
\medskip

\medskip
When $u$ is a {\it subpolynomial} of $p$, Lemma 1 tells us that the operation $X$ maps 
maps ${\cal H}$ to itself; it need not preserve degree of
course. The operation defined by replacing $p$ with $X(p)$ 
is a simple special case of a tensor product
operation defined in [D1]. 

\medskip
{\bf Definition 2}. An element $p$ of ${\cal H}(n,d)$ is called a 
{\it generalized Whitney mapping} if there exist elements $g_0,..., g_d$ of ${\cal H}(n)$ such that

1) $g_0 =1$ and $g_d = p$.

2) For each $j$, the degree of $g_j$ is $j$.

3) For each $j > 0$, we have $g_j = X(g_{j-1})$.

We say that $g_0,..., g_d$ defines a {\it Whitney Chain} from $1$ to $p$.

\medskip
At each step along the way of a Whitney chain, we replace $g_j$ with $g_j-u +su$, where
$u$ has degree $j$, and hence $g_k$ has degree $k$ for all $k$.

\medskip

{\bf Example 1}. The polynomial $x+xy+xy^2 + y^3$ is a generalized Whitney mapping with $d=3$.
 We have

$$ g_0 = 1 \mapsto g_1 = x+y \mapsto x + y(x+y) = g_2 =
x+xy + y^2 $$
$$ \mapsto x+xy + y^2(x+y) = g_3 = p = x+xy+xy^2 + y^3. \eqno (11) $$

We can rewrite (11) using the operation $X$:

$$ x+xy+xy^2 + y^3 = X(x+xy +y^2) = X(X(x+y))= X(X(X(1))).$$

\medskip
{\bf Lemma 2}. Suppose that $p \in {\cal H}(n,d)$ is a generalized Whitney mapping.
Then $N(p) \ge d(n-1) + 1$.

Proof. We induct on $d$. When $d=0$ we have $p=1$ and the conclusion holds. Suppose that
we know the result in degree $d-1$. Then $p = X(g) = g-u + su$, where $g$ is of degree $d-1$.
By the induction hypothesis, $N(g) \ge (d-1)(n-1) +1$. Suppose first
that $u$ consists of a single monomial $m$. Then $m$ is eliminated in passing from $g$
to $g-u$, but $m$ gets replaced with the $n$ new monomials $x_1 m,...,x_n m$.
Thus 
$$ N(X(g)) \ge N(g) + n-1 \ge (d-1)(n-1) + 1+ n-1 = d(n-1) + 1. \eqno (12) $$
If $u$ consists of several monomials, then because the coefficients
are nonnegative (12) remains true. $\spadesuit$

\medskip

We make a few simple remarks. First, the operation in (10) can be generalized
by replacing $s$ with any element of ${\cal J}$. Next, we show below that not all elements
of ${\cal H}(n)$ are generalized Whitney maps. On the other hand, if we allow negative coefficients
along the way, all such maps can be built up in this way.  We provide a simple
example.

\medskip
{\bf Example 2}. Consider $p(x,y) = x^3 + 3xy + y^3$.
Then $p \in {\cal H}(2,3)$. We can write $p=X_3 (X_2 (X_1(1)))$ as follows:

$$ 1 \mapsto s \mapsto s^2 = 3xy + s^2 -3xy  \mapsto 3xy +(s^2-3xy)s = 3xy+x^3+y^3 =
p(x,y).$$

In the notation (10), we have $g=s^2$ and $u=s^2 - 3xy$.
In using $s^2 -3xy$ we introduced a negative coefficient which was eliminated by the last
multiplication by $s$. One can easily show that we cannot construct $p$ by 
iterating this process while keeping all coefficients nonnegative.
As we stated above, if we allow negative coefficients along the way, 
then all elements of ${\cal H}(n)$ are obtained via iterations analogous to those in Example 2.
We now prove this assertion. 

Proposition 1 describes all elements of ${\cal H}(n)$ via
{\it undoing} the operation in (10). Proposition 2 uses only the operation (10)
but requires negative coefficients at intermediate steps. In proving these results
it is convenient to expand polynomials
in terms of their homogeneous parts. When $p$ is of degree $d$ we write
$$ p = \sum_{j=0}^d p_j, \eqno (13) $$
where each $p_j$ is homogeneous of degree $j$, and
we allow the possibility that $p_j=0$. 

\medskip

{\bf Proposition 1}. Suppose $p \in {\cal H}(n,d).$  Then there is an integer $k$ such
that 

$$ s^d = X^k(p) = \sum_{j=0}^d p_j s^{d-j}. \eqno (14) $$

Proof. Write $p=\sum p_j$ as in (13).  Suppose first that $p$ is not already
homogeneous. It is evident for each $j$ that
$p_j \subset p$. Let $\nu$ be the smallest index for which $p_\nu \ne 0$.
Then $p_\nu$ is a subpolynomial of $p$ and we may consider $X(p)$ defined as in (10) by

$$ X(p) = (p-p_\nu) + s p_\nu. $$
Then $X(p)$ also lies in ${\cal H}(n,d)$, and $X(p)$ vanishes to higher order than $p$
does. We iterate Lemma 1 in this way until we obtain the polynomial

$$ h = \sum_{j=\nu}^d s^{d-j} p_j, \eqno (15) $$
which lies in ${\cal H}(n,d)$. Now $h$ is homogeneous of degree $d$. 
The only homogeneous polynomial of degree $d$ that is identically  
equal to unity on the hyperplane $\{x: s(x) = 1\}$ is $s^d$. Therefore (14) holds.
$\spadesuit$.
\medskip
Formula (14) holds even when $p \in {\cal J}$, and we obtain the following
version where negative coefficients are allowed.

\medskip
{\bf Proposition 2}. Suppose $p \in {\cal J}(n,d)$. Then there is a finite list of maps
$X_1$,...$X_t$ from ${\cal J}$ to itself, of the form (17), such that 
$$ p = X_t \circ X_{t-1} \circ ...X_1(1). \eqno (16) $$

$$ X_j (v) = (v-r) + sr \eqno (17) $$

Proof. We induct on the degree. When the degree is zero, the only example
is $p=1$. Suppose that the result holds for all
elements of ${\cal J}(n,k)$ for $k \le d-1$. Let $p \in {\cal J}(n,d)$.
We expand $p$ into its homogeneous parts as above, and use (14) 
to rewrite the highest order part $p_d$. We obtain for a homogeneous
polynomial $r$ of degree $d-1$ that

$$ p = \sum _{j=0}^{d-1} p_j + p_d = \sum _{j=0}^{d-1} p_j + s^d - 
\sum _{j=0}^{d-1} p_j s^{d-j} = $$
$$ \sum _{j=0}^{d-1} p_j + s(s^{d-1} - \sum_{j=0}^{d-1} p_j s^{d-j-1}) = 
\sum _{j=0}^{d-1} p_j +sr = (p-p_d) +sr. \eqno (18) $$

Note that $p - p_d +r \in {\cal J}(n,d-1)$ and hence by the induction hypothesis
it can be factored as in (16). Since 
$$ p = (p-p_d) + sr = X (p - p_d + r), \eqno (19) $$
the induction step is complete. $\spadesuit$

\medskip
We repeat one subtle point regarding Proposition 2. Given $p \in {\cal H}(n,d)$,
it follows from (19) that
there exists $r$ of degree $d-1$ such that $p=u+sr$. In general
neither $r$ nor $u$ must have nonnegative coefficients.
The next mapping provides both an example where negative coefficients arise and an
example where the sharp bound from Theorem 0 arises without group invariance.
\medskip
{\bf Example 3}. Put $p(x,y) = x^7 + y^7 + {7 \over 2}x^5 y + {7 \over 2}x y^5 + {7
\over 2} x y$. Then $p \in {\cal H}(2,7)$. Following the proof of Proposition 2
we obtain

$$ p(x,y) = p_2(x,y) + p_6(x,y) + p_7(x,y) = $$
$$ p_2(x,y) + p_6(x,y) + (x+y)^7 - (x+y)^5
p_2(x,y) - (x+y)p_6(x,y),$$
and hence

$$ p = p_2 + p_6 + s (s^6 - p_2s^4 - p_6) = p - p_7 + sr. \eqno (20)$$
Here $r=s^6 - p_2s^4 - p_6$. Expanding $r$ yields

$$r(x,y) = x^6 - x^5 y + x^4 y^2 - x^3 y^3 + x^2 y^4 - x y^5 + y^6, \eqno (21) $$
which has negative coefficients.  Furthermore, $(p-p_7) + r$ has a negative coefficient.

\medskip
The operation $X$ replaces $u$ with $u-r + sr$. When we want to remind the reader
that we want both $r$ and $u-r$ to have nonnegative coefficients, we write $W$ instead of $X$.
To repeat, we cannot realize all 
elements of ${\cal H}(n)$ by successive application of $W$. 
We write ${\cal W}$ for the subset of ${\cal H}$ that can be obtained by
repeated application of the operation $W$ beginning with the constant function $1$.
We give one more simple example. Let $n=3$ with variables $(x,y,z)$. Applying
$W$ always to the ``last'' monomial, we obtain:

$$ W^3(1) = W^2(x+y+z) =W(x+y+xz+yz +z^2) = x+y+xy+xz+xz^2+yz^2 +z^3. \eqno (22) $$

We next give, without proof,
another  example of an element of ${\cal H}(n)$ that is not in ${\cal W}$.
The polynomial defined by (23) occurs also in Example 5. It plays an important role
because it satisfies the sharp estimate from Problem 1, yet it is not in ${\cal W}$.
In some sense it can exist because the dimension $3$ is too small
for stabilization to have taken place.

$$ x^3 + 3xy + 3 xz + y^3 + 3 y^2 z + 3 y z^2 + z^3. \eqno (23) $$
Observe that both (22) and (23) are of degree $3$, and each has $7$ monomials.

It is easy to see that polynomials formed by the process in (22) have
$N=d(n-1)+1$ terms. The first author has conjectured, for $n\ge 3$, that the inequality

$$ N \ge d(n-1) + 1 \eqno (24) $$
always holds. Theorem 2 yields
this inequality for all $n$ that are large enough relative to $d$. 
Given $d$, for such sufficiently large $n$ we prove a
stronger result by identifying all polynomials for which equality holds in (24);
these are precisely the generalized Whitney polynomials.
The stronger assertion fails in dimension three, but we believe that (24) still holds.

We next observe that there are always at least $n$ terms of degree $d$.

\medskip
{\bf Lemma 3}. Suppose $f \in {\bf R}[x_1,...,x_n]$ and $f$ is not identically $0$.
Then the polynomial $sf$ has at least $n$ monomials.

Proof. We claim first it suffices to assume that $f$ is homogeneous.
Assuming that the homogeneous case is known, then write
$ f = f' + f_d $, where $f_d$ consists of the highest degree terms.
Then $sf = sf' + s f_d$, where $s f_d$ has at least $n$ terms. 
All the terms in $sf'$ are of lower degree 
and hence cannot cannot cancel the terms in $s f_d$.  Thus the claim holds.

To prove the homogeneous case we proceed by induction on $n$.
When $n=1$ the result is trivial.
Suppose $n\ge 2$ and the result is known in $n-1$ variables. Given a homogeneous
$f$ in $n$ variables we write

$$ f(x) = x_n^d  f( { x_1 \over x_n},..., { x_{n-1} \over x_n}, 1) = 
 x_n^d  f(y_1,...,y_{n-1},1). \eqno (25) $$
It follows that

$$ s(x) f(x) = (y_1 + ... + y_{n-1} + 1) x_n^{d+1} f(y_1,...,y_{n-1}, 1). \eqno (26) $$
The number of terms in $sf$ is the same as the number of terms in the right-hand
side of (26) after dividing by $x_n ^{d+1}$.
Hence the number of terms in $sf$ is the number of terms in

$$ (y_1 + ...+ y_{n-1}) f(y_1,...,y_{n-1}, 1)  + f(y_1,...,y_{n-1}, 1). \eqno (27) $$
The first expression in (27) has at least $n-1$ terms by the induction hypothesis
and the second expression has at least one additional term.
$\spadesuit$.

\medskip

{\bf Corollary 1}. If $d>0$ and $p \in {\cal J}(n)$ has degree $d$, then 
$p$ has at least $n$ terms of degree $d$.

Proof. We write $ p = p' + p_d = p' + sr_{d-1}$ by (19). By Lemma 3, 
$s r_{d-1}$ has at least $n$ terms of degree $d$. $\spadesuit$
\medskip
We will close this section by proving Proposition 4 below. 
First we introduce a Veronese mapping
$\phi_{n-1} :{\bf R}^2 \to {\bf R}^n$ defined by

$$ \phi_{n-1}(u,v) = \left( u^{n-1},..., {n-1 \choose j}u^j v^{n-1-j},..., v^{n-1} \right). \eqno (28) $$

The Binomial Theorem shows that
the sum of the components of $\phi_{n-1}$ is $(u+v)^{n-1}$.
Therefore $\phi_{n-1}$ maps the line given by $u+v=1$ to the hyperplane $H$.

Let $p:{\bf R}^n \to {\bf R}$ be a function.
The pullback $\phi_{n-1}^*(p)$ is the composite function defined on ${\bf R}^2$
by $(u,v) \to p(\phi_{n-1}(u,v))$. We easily obtain the following
simple facts. 

\medskip
{\bf Proposition 3}. If $p \in {\cal H}(n,d)$, then 
$\phi_{n-1}^*( p) \in {\cal H}(2,d(n-1))$.
Furthermore $N(\phi_{n-1}^*(p)) \le N(p)$.

Proof. That $\phi_{n-1}^*(p)$ has degree $(n-1)d$ follows because $\phi_{n-1}$ is
homogeneous and the positivity of all coefficients prevents cancellation.
By the comment after (28) $$ \phi_{n-1}^* (s)(u,v) = s(\phi_{n-1}(u,v)) = (u+v)^{n-1}, $$
and thus $\phi_{n-1}$ maps the line given by $u+v=1$ to the hyperplane $H$.
Since $p=1$ on $H$, we see that $\phi_{n-1}^*(p) = 1$ on $u+v=1$.
Since all the coefficients
are all nonnegative, $\phi_{n-1}^*( p) \in {\cal H}(2,d(n-1))$.
Finally, we cannot increase the number of terms by a monomial substitution, and hence
$N(\phi_{n-1}^*(p)) \le N(p)$.  $\spadesuit$

\medskip

The proof of Proposition 3 
uses the nonnegativity of the coefficients.
For example, the pullback of the polynomial $x_2^2 - 4x_1x_3$ to 
$(u^2,2uv,v^2)$ vanishes. Without assuming nonnegativity of the coefficients
we cannot therefore conclude that the degree of $\phi_{n-1}^* (p)$ is $(n-1)d$.
The same example shows that pulling back via $\phi_{n-1}$ can
decrease the number of terms.

\medskip
{\bf Propsition 4}. Suppose $p \in {\cal H}(n,d)$. Then

$$ d(p) \le {2N(p)-3 \over n-1}. \eqno (29) $$

Proof. By Proposition 3 and Theorem 0 we obtain the chain of inequalities:

$$ d(p) = {d(\phi_{n-1}^*(p)) \over n-1} \le { 2N(\phi_{n-1}^*(p))-3 \over n-1}
\le {2N(p) - 3 \over n-1},$$
which gives the desired conclusion. $\spadesuit$

\medskip
The inequality in Proposition 4 is not sharp unless $n=2$.
When $n\ge 3$ the bound (5) obtained in Theorem 1 is smaller 
than the right-hand side of (29). For a given polynomial we can sometimes
obtain a better bound by pulling back via a mapping other than the Veronese.
We illustrate with a simple example. Define the mapping $p \in {\cal H}(3,7)$ by

$$ p(x,y,z) = x^3 + 3x(y+z) + (y+z)^3. $$
We have $d(p)=3$ and $N(p)=7$.
Pulling back via the Veronese mapping $\phi$ given by $\phi(u,v)=(u^2,2uv,v^2)$ gives
an element of ${\cal H}(2,6)$ with $7$ terms. The inequality 
$$ d(\phi^*(p))=6 \le 11 = 2N(\phi^*(p))-3 $$
is not sharp. Pulling back via the mapping given by $\psi(u,v) = (u^3, {\sqrt 3}uv,v^3)$ yields
an element of ${\cal H}(2,9)$ with $6$ terms, and therefore we obtain the sharp result
$$ d(\psi^*(p))= 9 = 2N(\psi^*(p))-3. $$
This discussion motivates the technique used to prove Theorem 1.

\medskip
{\bf IV. Optimal polynomials}
\medskip

We call an element $p$ of ${\cal H}(n,d)$ {\it optimal} if, for every $f \in 
{\cal H}(n,d)$, we have $N(f) \ge N(p)$. By Theorem 0, for $d$ odd, 
$p \in {\cal H}(2,d)$ is optimal if and only if $d=2N(p)-3$. 
The polynomials in (7) are optimal.
We hope to prove when $n\ge 3$ that $p \in {\cal H}(n)$ is optimal
if $N(p)=(n-1)d(p) + 1$. We can easily exhibit polynomials in
${\cal H}(n,d)$ for $n\ge 3$ satisfying this equality.

\medskip
{\bf Example 4}. Let $ s'(x)=\sum_{j=1}^{n-1} x_j$. We define $g_d$ by

$$ g_d(x) = x_n^d + s'(x) \sum_{k=0}^{d-1} x_n^k. \eqno (30) $$

It is evident from (30) and the finite geometric series
that $g_d \in {\cal W}$ and $N(g_d)=(n-1)d + 1$.

\medskip
{\bf Remark}. For a given $n$ and $d$ there are only finitely many optimal examples,
but typically there is more than one.  When $n=2$, for example, the first author has shown
the following fact. There are infinitely many $d$ for which there
exist optimal examples other than those given in (7) and those obtained by interchanging
the roles of $x$ and $y$. We omit the proof here.
Example 3 gives such an optimal polynomial of degree $7$.
\medskip
As mentioned above it is possible to improve Proposition 4 by pulling back
to the optimal examples in two dimensions.
We illustrate by establishing the next two Lemmas.

\medskip
{\bf Lemma 4}. Suppose $n\ge 2$ and $p \in {\cal H}(n,d)$. 
If $p$ contains a monomial in one or two variables of degree $d$, then 

$$ d \le { 2N-3 \over 2n-3}. \eqno (31) $$

Proof. After renumbering we may assume that $p$ contains either $x_1^d$
or $x_1^a x_2^b$ where $a+b=d$. Set $D=2n-3$. We pull back
using the optimal map $\phi$ induced by $p_D$ as defined in (7). 
Order the variables such that $x_1 =u^D$ and $x_2 = v^D$.
In either case we are guaranteed a term in $\phi^*(p)$ of degree $Dd$.
Following reasoning similar to the proof of Proposition 4 we obtain

 $$ d(p) = {d(\phi^*(p)) \over D} \le { 2N(\phi ^*(p))-3 \over D}
\le {2N(p) - 3 \over D} = { 2N-3 \over 2n-3}, \eqno (32) $$
which gives (31). $\spadesuit$.
\medskip

By assuming that the highest degree part of $p$ contains monomials
involving few of the variables we can generalize the preceding proof.
We give two of several possible versions.

\medskip

{\bf Lemma 5}. Suppose $n\ge 2$ and $p \in {\cal H}(n,d)$. If $p$
contains the monomial $m =x_1^{a_1}...x_k^{a_k}$ of degree $d$, where $k\ge 2$,
then  the following hold:

$$ d(p) \le { 2N-2k+1 \over 2n-2k+1}. \eqno (33) $$

$$ d(p) \le { 2N-3 + \sum_{j=3}^k (j-2) a_j \over 2n-3}. \eqno (34) $$

Proof. First we prove (33).
We set $x_j = {\lambda \over k-1}$ for $2 \le j \le k$.
In doing so we replace $k-1$ terms with one term,
thus killing $k-2$ terms. We also decrease the number of variables by $k-2$.
We now pullback as in the proof of Lemma 4 (or use Lemma 4 directly)
to see that 

$$  d(p) \le { 2(N-(k-2))-3 \over 2(n-(k-2))-3} = { 2N-2k+1 \over 2n-2k+1}.$$
We have proved (33).

The proof of (34) also involves pulling back to the optimal polynomials
in two dimensions. We first set $D=2n-3$, and consider the mapping
$\phi$ induced by $p_D$ as defined in (7), where the coordinates are ordered such that

$$ (x_1,x_2,x_3, x_4,...) = (u^D, v^D, c_1 u^{D-2}v, c_2 u^{D-4} v^2,... ) =
\phi(u,v). $$
Pulling back the monomial $m$ then guarantees a term of degree

$$ a_1D + a_2 D + a_3(D-1) +... + a_k(D-k+2) = D \sum_{j=1}^k a_j - \sum_{j=3}^k (j-2)a_j$$
in $\phi^*(p)$. Since the sum of the $a_j$ is $d$ we obtain

$$ dD-\sum_{j=3}^k (j-2)a_j \le d(\phi^*(p)) \le  2N(\phi^*(p)) - 3 \le 2N(p)-3, \eqno (35) $$
and hence

$$  d(p) = d \le {2N(p)-3 + \sum_{j=3}^k (j-2)a_j \over D} = 
{2N(p)-3 + \sum_{j=3}^k (j-2)a_j \over 2n-3}. \eqno (36) $$ 
Thus we have proved (34).  $\spadesuit$

\medskip
The proof of (34) when $k=2$ is essentially the same as the proof of Lemma 4.
The proof of (34) gives the strongest result by taking $D$ as large as possible;
$D=2n-3$ is the largest number for which $\phi$ takes values in $n$-space, a
requirement for the proof to make sense. Thus the
choice of $D$ itself relies on Theorem 0.

Let us write $E= \sum_{j=3}^k (j-2)a_j$. Our next result provides a general bound
for $d(p)$ in terms of $N(p)$ in all cases. We do so by estimating 
the {\it excess} $E$ in terms of $d$ and $n$. From Theorem 1 we obtain the weaker
asymptotic bound

$$ d(p) \le {4 \over 3} {2N(p) - 3 \over 2n-3} $$ 
as $n\to \infty$. Our main result, Theorem 2, provides the 
sharp asymptotic result $d \le {N-1 \over n-1}$ when $n$ is large relative to $d$.
On the other hand Theorem 1 holds for all $n$ and its proof is much simpler,
but it is sharp only in two dimensions.

\medskip {\bf Theorem 1}. Suppose $p \in {\cal H}(n,d)$. Then

$$ d(p) \le {2n (2N(p)-3) \over 3n^2 -3n -2} \le {4 \over 3} {2N(p)-3 \over 2n-3}. \eqno (37) $$

Proof. We begin with the estimate
$$  d(p) \le {2N(p)-3 + \sum_{j=3}^k (j-2)a_j \over 2n-3}\eqno (38) $$ 
from Lemma 5. For notational ease we rewrite (36) as

$$ d(p) \le F + {E \over D} \eqno (39)$$
where $F = {2N-3 \over 2n-3}$.
We may assume $k\ge 2$ and that $a_1 \ge a_2 \ge ...\ge a_k$. We obtain 

$$ {E \over D}= { \sum_{j=3}^k (j-2)a_j \over D} \le {d \over Dk} \sum_{j=3}^k (j-2) = 
{d \over Dk} {k-1 \choose 2}. \eqno (40) $$
Since $k\le n$, we obtain from (40) the upper estimate

$$ {E \over D} \le {d \over nD} {n-1 \choose 2} = c(n) d, \eqno (41)$$
where the expression $c(n)$ is defined by
$$ c(n) = {{n-1 \choose 2} \over n(2n-3)}. \eqno (42) $$
One easily shows that $c(n) < 1$. Therefore 
(39) yields

$$ d(p) \le F + {E \over D} \le F + c(n) d(p) $$
and hence
$$ d(p) \le {1 \over 1-c(n)} F = 
{2N-3 \over 2n-3} {1 \over 1-c(n)} = {2n (2N(p)-3) \over 3n^2 -3n -2}. \eqno (43) $$
We have bounded $d$ in terms of $N$ and $n$. It is elementary to verify
for $n\ge 2$ that 
$$ {2n \over 3n^2-3n-2} \le {4 \over 3(2n-3)}, $$
and therefore the inequality on the far right-hand side of (37) holds.  $\spadesuit$
\medskip

We pause to mention an explicit optimal example.

$$ p(x,y,z) = x + y + z^2 + xz + y^2z + yz^2 + xyz(x+y+z). \eqno (44) $$

The polynomial in (44) is of degree $4$, but each term of degree $4$ involves
all three of the variables and thus Lemma 4 is not useful.
Note that $N(p)=9$. By Proposition 5, 
nine is the smallest possible number of terms for an element in ${\cal H}(3,4)$. 

\medskip
Before turning to Proposition 5, which is proved below and verifies the conjecture 
(3) from Problem 1 for degree up to 4, we
briefly discuss one parameter families of mappings.  The following proposition 
will be proved and developed in [L].  A one-parameter family of polynomials is defined by 
$$ p_\lambda(x) = \sum c_\alpha (\lambda) x^\alpha, \eqno (45) $$
where each map $\lambda \to c_\alpha(\lambda)$ 
is a continuous function of a real parameter $\lambda$. 
One simple example of a one-parameter family is given by the convex combination
$f_\lambda = \lambda p + (1-\lambda)q$ of elements $p$ and $q$ of ${\cal H}(n,d)$.
We observed earlier that $f_\lambda \in {\cal H}(n,d)$ as well.

\medskip
{\bf Proposition L}. Let $p_t$ denote a one-parameter family of elements of 
${\cal H}(n,d)$. Suppose that $N(p_t)$ is constant for $t$ in an open interval.
Then $p_t$ is optimal for no $t$.

\medskip

We next include some information
which supports the conjectured sharp bound.
The proofs of the four statements in the following result become increasingly
elaborate as the codimension increases. We therefore provide detailed proofs
of statements 0), 1), and 2) but only an outline of the proof of 3).
The proofs of 0) and 1) are easy; the proofs of 2) and 3)
first use combinatorial reasoning to make Lemma 4 applicable and then use
additional combinatorial reasoning to improve the bound
from Lemma 4 in these special cases. The bounds in this result are interesting in the
context of CR mappings between spheres.

\medskip
{\bf Proposition 5}. Suppose $p \in {\cal H}(n,d)$ for $n\ge 3$. Then

0) If $N(p) < n$, then $d = 0$.

1) If $N(p) < 2n-1$, then $d \le 1$.

2) If $N(p) < 3n-2$, then $d \le 2$.

3) If $N(p) < 4n-3$, then $d \le 3$.

Proof: The contrapositive of 0) is easy. When $d \ge 1$ there must 
be at least $n$ distinct monomials of degree $d$, by Corollary 1.

We call terms of the form $x_i^k$ {\it pure} terms, and we call monomials
depending on at least $2$ variables {\it mixed} terms. 
By pulling back to the one-dimensional case in $n$ ways
(by setting $n-1$ of the variables equal to zero),
we note that there must be at least $n$ distinct pure terms. 
If $d=1$ then all the terms are pure terms and $p=s$. We may therefore assume that
$d\ge 2$ in proving the rest of the statements. 

The proof of 1) proceeds as follows.
If no pure term is of degree at least $2$ then as above $p=s$. We may thus assume that 
the monomial $x_1^a$ occurs for some $a \ge 2$. By setting all variables 
except $x_1$ and $x_j$ equal to $0$, we see that a mixed monomial $x_1^k 
x_j^l$ must occur for $2\le j\le n$. Hence we have at least $n-1$ mixed terms.
Counting also the $n$ pure terms shows that $N(p) \ge (n-1) + n$ and we obtain 1).

If $d=2$ then 2) holds. We therefore assume $d \ge 3$
when proving 2). We must then show that $N \ge 3n-2$. There are
two cases: 

If $x_1^a$ is the only pure term of degree greater than $1$ then $p$ must be 
equal to $x_1 r(x) + s - x_1$, for some $r(x) \in {\cal H}(n,d)$. The 
polynomial $r$ has $n-1$ fewer terms than $p$ does and it must have degree at least $2$.
Applying 1) shows that 
$N(r) \ge 2n-1$ and hence $N(p) \ge (2n-1) +(n-1) = 3n-2$. 
Thus 2) holds in this case.
 
The remaining case of 2) is when at least two pure terms of degree at least $2$ occur.
Hence we assume that $x_2^b$ occurs as well, with $b \ge 2$.
We then have at least $2(n-2)+1$ mixed terms and $n$ pure terms for a 
total of $3n-3$. We want $N\ge 3n-2$. Let us therefore 
assume for the purpose of contradiction that 
there are no other terms. For $d\ge 3$ the only element of ${\cal H}(2,d)$
that has at most 3 distinct monomials is $u^3 + 3uv + v^3$. Hence  
all pure terms must be of degree $3$ and we obtain
$$ p(x) = \sum_{j=1}^n x_j^3 + 3 \sum_{i \neq j} x_jx_i. \eqno (46) $$
We claim that the polynomial in (46) is not in ${\cal H}(n,3)$ unless $n=2$. 
To verify the claim we note that $p({1 \over n},...,{1 \over n}) > 1$ when $n \ge 3$. 
Thus 2) holds in this case, and hence in general.

To prove 3) we assume $N \le 4n-4$.
If Lemma 4 does not apply, then there is no term of degree $d$
involving at most two of the variables. We must then have at least
$n$ terms of top degree, $n$ additional pure terms, and (as above) at least $2n-3$
additional mixed terms involving two variables. The total is $4n-3$ and thus
$N \ge 4n-3$. We may therefore assume Lemma 4 applies. In particular $d \le 4$.

We proceed by contradiction. Assume $d=4$.
We consider the cases $N\le 4n-5$ and $N=4n-4$ separately.
If $N \le 4n-5$ we obtain a contradiction as follows:
By Lemma 4, 

$$ d(2n-3) +3 \le 2N. $$
Including the information on $N$ and $d$ yields

$$ 4(2n-3) + 3 \le d(2n-3) + 3 \le 2N \le 2(4n-5) $$
from which we obtain the contradiction $-9 \le -10$. 
Thus, for $N\le 4n-5$ we have $d \le 3$.

The remaining case is when $N(p) = 4n-4$ and $d=4$.
There are two subcases.
First suppose that $n \ge 4$. 
As argued above we can assume that 
there exist pure monomials in $x_1$ and $x_2$ of degree greater than 1.
Setting in turn $x_1 = 0$ and $x_2 = 0$
we get polynomials in $n-1$ variables with at least $n$ fewer terms.  Thus
these polynomials must have degree at most 3.
The
top degree terms must be divisible by $x_1 x_2$, and thus $p_4 = s(x) x_1 x_2 q(x)$,
where $q$ is homogeneous of degree 1.  We can easily check that $q$ 
must have all positive coefficients, and we can undo an operation $X$ to 
reduce to a previous case.

The other subcase is when $n=3$, $N(p) = 4n-4=8$ and $d = 4$. We claim that 
no polynomial in ${\cal H}(3,4)$ has exactly $8$ distinct monomials. There are only 
finitely many possibilities that need to be checked and we outline 
how to do this by hand.

If all terms of degree $4$ depend on $3$ variables, we undo and reduce to
a previous case to get a contradiction.  After 
renaming variables, we consider the polynomials
$p(x_1,x_2,0)$, $p(x_1,0,x_2)$, and $p(0,x_2,x_3)$. A counting argument shows
that the first two of these must
have exactly 4 terms and be of degree $4$, whereas the third must have 3
terms and must be of degree $3$ or less.  By a study of the 2-dimensional case
we see that $x_1^4$ must appear.  One can then check by hand
that the only possible configuration of degree $4$ terms is $x_1^3 (x_1 + x_2 +
x_3),$ and reducing to a previous case produces a contradiction.  $\spadesuit$

\medskip

\medskip
The following corollary supports the conjectured sharp bound for degree at most $4$.
We believe that these bounds are sharp for all degrees when $n\ge 3$. 
In the next section we establish this result when $n$ is large enough compared with $d$.
\medskip

{\bf Corollary 2}.  Suppose $n\ge 3$ and
$p \in {\cal H}(n,d)$. If $d \le 4$ or $N(p) < 4n-3$, then the following two estimates hold:
$$ N(p) \ge d(n-1) + 1 , $$
$$ d \le {N(p)-1 \over n-1} . \eqno (47) $$

\medskip

{\bf V. Whitney Mappings and the Proof of Theorem 2}.
\medskip

In this section we give conditions under which
a polynomial $p \in {\cal H}(n,d)$ in fact lies in ${\cal W}$.  
By Lemma 2 if $p \in {\cal W} \cap {\cal H}(n,d)$ then 
the desired bound $N(p) \ge d(n-1)+1$ holds.

The following theorem is the main result of this paper. It solves Problem 1 
when the domain dimension is large enough.

\medskip
{\bf Theorem 2}.
Fix $d$ and assume $n \ge 2d^2 + 2d$.  If $p \in {\cal H}(n,d)$
then $N(p) \geq (n-1)d + 1$.  Furthermore, if equality holds then
$p \in {\cal W}$.
\medskip

Before we prove Theorem 2
we give a simple condition guaranteeing that $p \in {\cal W}$. 
Let $x = (x',x_n) \in {\bf R}^{n-1} \times {\bf R}$ and define
$s'(x') := \sum_{j=1}^{n-1} x_j$.  We will say that $p$ is {\it affine in}
$x_n$ if we can write
$p(x^\prime, x_n) = a(x^\prime) + x_n b(x^\prime)$ for some polynomials $a$ and
$b$.

\medskip
{\bf Lemma 6}. 
If $p \in {\cal H}(n,d)$ and suppose $p$ is affine in $x_n$,
then $p \in {\cal W}$.

Proof. We induct on the degree $d$. When $d=1$ the result is obvious.
Suppose $d \ge 2$ and
that the result is known for such affine polynomials of degree $d-1$.
Assume $p(x^\prime, x_n) = a(x^\prime) + x_n b(x^\prime)$. By (18) from Proposition 2
we write $p= (p- p_d) + s r_{d-1}$. Equating the highest part of these expressions 
for $p$ gives

$$ a_d(x') + x_n b_{d-1}(x') =
\left( \sum_{j=1}^{n-1} x_j + x_n \right) r_{d-1}(x')
= s'(x') r_{d-1}(x') + x_n r_{d-1}(x'). \eqno (48) $$
Hence $r_{d-1} = b_{d-1}$ and $a_d= s' r_{d-1}$. Therefore

$$ p = p-p_d + s b_{d-1} = X(p-p_d + b_{d-1}). \eqno (49) $$

Note that $p-p_d + b_{d-1} \in {\cal H}(n,d-1)$. It is also affine in $x_n$ and hence
lies in ${\cal W}$ by the induction hypothesis. Thus $p \in {\cal W}$ as well. 
$\spadesuit$
\medskip

We now prove two simple results that we use in the proof of Theorem 2.
The reader should look back at Examples 1 and 4.

\medskip
{\bf Lemma 7}. Let $p \in {\cal H}(2,d)$ and
suppose that $p(x,y) = a(x) + y b(x)$.
Then $N(p) \ge d+1$. The monomial $x^d$ must appear
and $x^j y$ must appear for each $j$ with 
$0 \le j \le d-1$.
Furthermore, $p$ has exactly $d+1$ distinct monomials if and only if
$$ p(x, y) = x^d + y(x^{d-1} + \cdots + x + 1).$$

Proof.
By Lemma 6 we know $p \in {\cal W}$, and 
the statement follows by induction on $d$.
$\spadesuit$
\medskip

For two monomials $m_1 = x_1^{\alpha_1} \cdots x_n^{\alpha_n}$ and $m_2 =
x_1^{\beta_1} \cdots x_n^{\beta_n}$ we define the distance between them
by
$$
\delta(m_1, m_2) := \sum_j |\alpha_j - \beta_j|. $$
For monomials of the same degree 
$\delta(m_1, m_2)$ must be even.

\medskip
{\bf Lemma 8}.
Let $p \in {\cal H}(3,d)$, and suppose that
$p(x_1, x_2, x_3) = a(x_1,x_2) + x_3 b(x_1,x_2)$.
If two monomials $m_1(x_1,x_2), m_2(x_1,x_2)$ of degree $d$ occur in
$p(x)$ with $\delta(m_1, m_2) \ge 4$, then $p$ has at least $d+1$ distinct
monomials that depend on $x_3$.

Proof. 
It follows from Lemma 6 that $p \in {\cal W}$, and from
Lemma 7 that $p$ must have at least one monomial of
every degree that depends on $x_3$.  Since $\delta(m_1, m_2) \ge 4$ there
must be at least $2$ monomials of maximal degree that depend on $x_3$,
which gives at least $d+1$ monomials. $\spadesuit$
\medskip
For the rest of this section we assume $n \geq 2d^2 + 2d$.
In particular $n \ge 3$.
Let $p \in {\cal H}(n,d)$ and let $N=N(p)$. We assume both
that $N \leq d(n-1) + 1$ and that $p$ is optimal.  We will show
that $p$ must be a generalized Whitney mapping and thereby prove Theorem 2.

Let $m_1$ and $m_2$ be distinct monomials that occur in $p$. 
The main idea of the proof is to show that $\delta(m_1, m_2)$ must be equal to $2$.

Let $k$ be the number of distinct variables that occur in either $m_1$ or $m_2$. 
Then $2 \le k \le 2d$.
After renaming the variables if necessary we may assume that $m_1$ and 
$m_2$ are independent of $x_j$ for $j \ge k+1$.

We define new polynomials in ${\cal H}(2,d)$ and ${\cal H}(3,d)$
$$ P_{j}(\xi,x_j) :=
p\Big(\underbrace{{\xi \over k},\ldots,{\xi \over k}}_{k \ {\rm times}},
0,\ldots,0,x_j,0,\ldots\Big) , $$

$$ P_{ij}(\xi, x_i,x_j) :=
p\Big(\underbrace{{\xi \over k},\ldots,{\xi \over k}}_{k\  {\rm times}},
0,\ldots,0,x_i,0,\ldots,0,x_j,0,\ldots\Big) . \eqno (50)$$

\medskip
{\bf Claim}.
The polynomial $P_j$ is affine in $x_j$ for each $j \in \{k+1, \ldots,n \}$.

Proof. Seeking a contradiction we assume $k+1 \le l \le n$,
that $P_j$ is not affine for $k+1 \le j \le l$, and that $P_j$ is affine
for $l+1 \le j \le n$.

If $P_j$ is affine in $x_j$ then by Lemma 6 we have
$$
P_{j}(\xi,x_j) = c_1 \xi^d + c_2 \xi^{d-1} x_j + \cdots + c_d \xi
x_j + c_{d+1} x_j + q(\xi), $$
where $q$ is a possibly zero polynomial in $\xi$ of degree $d-1$
or less. If $P_{j}$ is not affine in $x_j$ then there must be at
least $\lceil {d-3 \over 2} \rceil$ terms by Theorem 0.

We will proceed to find a lower estimate for the number of
monomials of $p$, and we must take care not to count the same
monomial twice. We first count the monomial $m$. For each
$P_{j}$ where $k+1 \leq j \leq l$ we have at least $\lceil
{d+3 \over 2} \rceil - 1$ extra monomials and for each $P_{j}$ for
$j > k$ we get at least $d$ extra monomials.

For $P_{ij}$ where $k+1 \leq i < j \leq l$ we know that there must
be at least one monomial that depends on $x_i$ as well as $x_j$
(keep $\xi$ constant to see this),
and thus we get least $(l-k)(l-k-1)/2$ more monomials that we have
not counted yet.

For the same reason we can count one extra monomial depending on
both $x_i$ and $x_j$ for each possible choice $k+1 \leq i \leq l <
j \leq n$ so we get $(l-k)(n-l)$ more monomials.

When we add the number of all these monomials we obtain

$$ N \geq 1 + (l-k)\left(\left\lceil {d+3 \over 2} \right\rceil - 1 +
{l-k-1 \over 2} + (n-l)\right) + (n-l)d . \eqno (51) $$

By our assumption $l \ge k+1$.  If
$$(l-k)\left(\left\lceil {d+3 \over 2}
\right\rceil - 1 + {l-k-1 \over 2} + (n-l)\right) > (l-1)d, \eqno (52) $$
then $p$ cannot be optimal.  This happens when
$$ (l-k)(d - l-k + 2n) - 2(l-1)d > 0 . \eqno (53)$$

Fixing $k, d$ and $n$ the expression in (53) is concave down in $l$ and thus must
achieve a minimum if $l = k+1$ or $l = n$.  We know $2 \leq
k \leq 2d$ and so get two bounds for $n$:
$$n > {4d^2+3d+1 \over 2} ,$$

$$n > 5d. \eqno (54)$$
Our assumption that $n \ge 2d^2 + 2d $ implies both bounds
(noting that $d \geq 2$). We have proved the Claim.
\medskip

Now suppose for the sake of contradiction that 
$\delta(m_1, m_2)$ is at least $4$. Write $m_1 = \prod_{i=1}^{k} x_i^{r_i}$ and $m_2 =
\prod_{i=1}^{k} x_i^{s_i}$. By renaming the variables again if necessary we assume that there exists
an integer $t$ such that for $i = 1, \ldots , t$
we have that $r_i \ge s_i$, and for $i = t+1, \ldots , k$
we have $r_i \le s_i$. It follows from the claim
that for $j =  k+1 , \ldots , n$ the polynomial $P_j$, 
as defined in Equation (50), must be affine in $x_j$.

Let
$$
P(y,z,x_{k+1} , \ldots ,x_n) := p\Big(
\underbrace{{y \over t},\ldots,{y \over t}}_{t \ {\rm times}},
\underbrace{{z \over k-t},\ldots,{z \over k-t}}_{k-t \ {\rm times}},
x_{k+1}, \ldots , x_n \Big) .
$$
It follows
that $P$ has two terms of highest degree $y^{r_1}z^{r_2}$
and $y^{s_1}z^{s_2}$ with $r_1 > s_1 + 1$ and $r_2 < s_2-1$.
Therefore for every $j \in \{k+1, \ldots , n\}$, the
polynomial $P(y,z, 0, \ldots , 0 ,x_j , 0, \ldots , 0)$ is
a polynomial in three variables that satisfies the conditions of
Lemma 8, and hence it has at least $d+1$ terms
that depend on $x_j$. Hence $P$ (and thus also $p$)
has at least $(d+1)(n-2d) = dn + n -2d^2 - 2d$ distinct monomials.
We assumed that $n
\ge 2d^2 + 2d$, so the polynomial cannot be optimal, which contradicts our assumption.
Thus $\delta(m_1, m_2) = 2$.

By Corollary 1 there are at least $n$ terms of
highest degree.  It follows that the terms of highest degree must equal
$c s \cdot m$ for some constant $c$ and some monomial $m$ of degree $d-1$.
Recall that $s$ denotes the sum of the variables.

Thus we can undo the operation $X$ to obtain a new
polynomial of degree $d-1$, with exactly $n-1$ terms fewer than
$p$.  The reason is that $p$ is optimal; undoing the operation $X$
must create a new term of degree $d-1$ (otherwise multiplying that term
by $s$ would get a polynomial with fewer terms than $p$).
This new polynomial of degree $d-1$ must again be optimal, because if there
existed a polynomial of degree $d-1$ with fewer terms, we could apply
operation $X$ to it and again and invalidate the optimality of $p$.

An inductive argument with respect to the degree shows that $p$
must be obtained by starting with $s$ and
repeatedly multiplying one of the highest degree terms with $s$, in other
words, $p \in {\cal W}$.  We have completed the proof of Theorem 2. $\spadesuit$

\medskip

\medskip
{\bf VI. CR Mappings between Spheres}
\medskip
The results of this paper are closely related to a basic question in CR Geometry.
Let $f$ be a rational mapping from complex Euclidean
space ${\bf C}^n$ to ${\bf C}^N$, and suppose $f$ maps the unit sphere $S^{2n-1}$ in its domain
to the unit sphere $S^{2N-1}$. Can we give any estimate for the degree of $f$ in terms of $n$ and $N$? The degree of a rational map $f={p \over q}$ is defined to be the maximum of the degrees of 
$p$ and $q$, when $f$ is reduced to lowest terms. It is easy to show in this context [D3] that the degree of $f$ equals the degree of $p$.

Many of the results mentioned below do not begin by assuming that $f$ is rational. 
Instead they assume that $f$ is
a proper mapping between balls, and they make some regularity assumptions at the boundary
in the positive codimension case.
By the work of Forstneric ([F1] and [F2]), a proper mapping between balls (with domain dimension
at least $2$), with sufficient differentiability at
the boundary, must be a rational mapping. We therefore assume rationality in this section.

We return to the basic question of degree.
As in this paper, when $n=1$ the answer is no. Assume next that $n\ge 2$.
As in Proposition 5 of this paper, when $N < n$ we can conclude by elementary considerations
that $f$ must be a constant. When $N=n \ge 2$, Pincuk [P] proved that $f$ must either be
a constant or a linear fractional transformation,
and hence of degree at most $1$. Faran [Fa1]
showed that we can draw the same conclusion when $n \le N \le 2n-2$. When $n=2$ and $N=2n-1 =3$, Faran [Fa2]
showed that, up to composition with automorphisms of the ball on both sides, the map must be a monomial
mapping of degree at most $3$. Thus the rational mapping is of degree at most $3$ in this case.
In particular Faran discovered the mapping $(u^3, {\sqrt 3}uv, v^3)$ which is of
maximum degree from the two-ball to the three-ball, and is group-invariant. In [D2], [D3], and [D5] the 
first author studied the group invariance aspects of CR mappings, discovered the maps (7), 
and observed many connections to other branches of mathematics. 

Huang and Ji have investigated ([H] and [HJ]) aspects of the basic  question.
They have established, for example, when $3 \le n \le N=2n-1$, that the degree of a (rational mapping (between
spheres) is at most $2$, and they have discovered various conditions somewhat analogous to our work here
for guaranteeing partial linearity. One striking aspect of their work is that they do not assume rationality
and their regularity assumptions are minimal. All these papers involve the low codimension case. 
Meylan's [M] result gives the bound $d \le {N(N-1) \over 2}$ in any codimension,
when the domain dimension $n$ is assumed to be two. The paper [HJX] includes the
following result. Let $f$ be a rational proper mapping between balls 
of degree $2$. If $f$ has {\it geometric degree} $1$, then $f$ is a generalized
Whitney map. 

The expository paper [D4] includes the relationship of this complexity issue 
to a complex variables analogue of Hilbert's 17th
Problem, and includes the following result. Given a rational mapping ${p \over q}:{\bf C}^n \to {\bf C}^N$
that maps the closed unit ball into the open unit ball, we can find an integer $K$
and another rational mapping ${g \over q}:{\bf C}^n \to {\bf C}^K$
(with the same denominator) such that the mapping $({p \over q}, {g \over q})$ maps $S^{2n-1}$ to $S^{2(N+K)-1}$.
We must be able to choose $K$ large enough. Even for quadratic mappings and $n=2$,
we must chose $K$ to be arbitrarily large. Thus by placing no restriction on the target dimension,
we can create arbitrarily complicated
rational mappings between spheres. In future work we will show how the bounds in this paper, which arise by
considering monomial rather than rational maps, can to some extent be extended to the rational case. 

The first author has conjectured that the degree of a rational mapping 
sending $S^{2n-1}$ to $S^{2N-1}$ is at most ${N-1 \over n-1}$ when $n\ge 3$, and it 
is at most $2N-3$ when $n=2$.
The results in this paper show how to obtain sharp results in the special but nontrivial
case where the map is a monomial.

\medskip
{\bf References}
\medskip

[D1] D'Angelo, John P., Proper polynomial mappings between balls, Duke Math J.
57(1988), 211-219.
\medskip

[D2] D'Angelo, John P.,
{\it Invariant holomorphic mappings}, 
Journal of Geometric Analysis Vol. 6, No. 2 (1996), 163-179.
\medskip

[D3] John P. D'Angelo, Several Complex Variables and the Geometry of Real
Hypersurfaces, CRC Press, Boca Raton, 1993.

[D4] John P. D'Angelo,
{\it Proper holomorphic mappings, positivity conditions, and isometric imbedding},
J. Korean Math. Society (2003), 1-30.

[D5] John P. D'Angelo, {\it Number-theoretic properties 
of certain CR mappings}, Journal of Geometric Analysis, Vol. 14, No. 2 (2004), 215-229.

[DKR] John P. D'Angelo, \v Simon Kos, and Emily Riehl,
{\it A sharp bound for the degree of proper monomial
mappings between balls}, Journal of Geometric Analysis, Volume 13, Number 4 (2003),
581-593. 

[Fa1] J. Faran, {\it Linearity of proper holomorphic mappings in the 
low codimension case}, J. Diff. Geom. 24 (1986), 15-17.

[Fa2] J. Faran, {\it Maps from the two-ball to the three-ball}, 
Inventiones Math., 68 (1982), 441-475.

[F1] Franc Forstneric, {\it Extending proper holomorphic mappings 
of positive codimension},
Inventiones Math. 95 (1989), 31-62.

[F2] Franc Forstneric, {\it Proper rational maps: A survey},  Pp 297-363 in
Several Complex Variables: Proceedings of the Mittag-Leffler
Institute, 1987-1988, Mathematical Notes 38, Princeton Univ.
Press, Princeton, 1993.

[H] X. Huang, {\it On a linearity problem for proper maps 
between balls in complex spaces
of different dimensions}, J. Diff. Geom. 51 (1999), no 1, 13-33.

[HJ] X. Huang  and S. Ji, 
{\it Mapping $B_n$ into $B_{2n-1}$}, Inventiones Math. 145 (2001), 219-250.

[HJX] X. Huang, S. Ji, and D. Xu,
{\it Several results for holomorphic mappings from ${\bf B}_n$ to 
${\bf B}_N$}, Contemporary Math. 368 (2005), 267-292.

[L] Ji{\v r}\'\i\ Lebl,
Singularities and Complexity in CR Geometry, PhD thesis, University of California, San Diego, 2007.

[M]  Francine Meylan, {\it Degree of a holomorphic map between unit balls from ${\bf C}^2$ to
${\bf C}^n$}, Proc. A.M.S., Vol. 134, No. 4, 1023-1030.

[P] S. I. Pincuk, {\it On the analytic continuation of holomorphic mappings}, 
Math USSR-Sb. 27 (1975), 375-392.

\end